\documentclass[opre, nonblindrev]{informs3notitle}
\OneAndAHalfSpacedXI 

\usepackage{amsmath}
\usepackage{amssymb}
\usepackage{subfigure}
\usepackage{bm}
\usepackage{url}
\usepackage{graphicx}
\usepackage{epstopdf}
\usepackage{epsfig,psfrag}

\usepackage{verbatim} 
\usepackage{subfigure}
\usepackage{algorithmic}
\usepackage{enumerate}
\usepackage{color}
\usepackage{pdfsync}
\usepackage{hyperref}

  \usepackage{natbib}
 \bibpunct[, ]{(}{)}{,}{a}{}{,}%
 \def\bibfont{\small}%
 \def\bibsep{\smallskipamount}%
 \def\bibhang{24pt}%
\global\long\def\zp{\mathbb{Z}_{+}}
\global\long\def\N{\mathbb{N}}

\global\long\def\pb{\mathbb{P}}
\global\long\def\calf{\mathcal{F}}

\global\long\def\E{\mathbb{E}}

\global\long\def\rp{\mathbb{R}_{+}}

\global\long\def\bproof{\emph{Proof.}~}
\global\long\def\arr{\mathcal{A}}
\global\long\def\ser{\mathcal{S}}
\global\long\def\sk#1#2{\stackrel{(#1)}{#2}}

\global\long\def\ch{c_h}
\global\long\def\cl{c_l}
\global\long\def\cnt{\mathcal{N}}
\global\long\def\hc{\mathcal{H}}
\global\long\def\jc{\mathcal{J}}

\global\long\def\eva{\mathcal{E}_1}
\global\long\def\evb{\mathcal{E}_2}
\global\long\def\evc{\mathcal{E}_3}
\global\long\def\evd{\mathcal{E}_4}
\global\long\def\eve{\mathcal{E}_5}

\global\long\def\qop{\mathcal{Q}^*(\lambda, \wlam)}
\global\long\def\qlam{q_\lambda}
\global\long\def\wlam{W_\lambda}

\global\long\def\lleq{\preccurlyeq}
\global\long\def\ggeq{\succcurlyeq}
\global\long\def\asgeq{\succeq}
\global\long\def\asleq{\preceq}

\global\long\def\bsep{\,\big|\,}

\global\long\def\nnb{\nonumber}

\usepackage{natbib}
 \bibpunct[, ]{(}{)}{,}{a}{}{,}%
 \def\bibfont{\small}%
 \def\bibsep{\smallskipamount}%
 \def\bibhang{24pt}%
\TheoremsNumberedThrough     

\EquationsNumberedThrough

\begin{document}


\RUNAUTHOR{Xu}

\RUNTITLE{Necessity of Future Information in Admission Control}

\TITLE{Necessity of Future Information in Admission Control}

\date{\today}


\ARTICLEAUTHORS{%
\AUTHOR{Kuang Xu}
\AFF{Stanford University \\
Graduate School of Business\\
Stanford, CA 94305}
\EMAIL{kuangxu@stanford.edu} 
} 

\ABSTRACT{We study the {necessity} of predictive information in a class of queueing admission control problems, where a system manager is allowed to divert incoming jobs up to a fixed rate, in order to minimize the queueing delay experienced by the admitted jobs.  

\cite{SSX12} show that the system's delay performance can be significantly improved by having access to {future information} in the form of a lookahead window, during which the times of future arrivals and services are revealed. They prove that, while delay under an optimal online policy diverges to infinity in the heavy-traffic regime, it can stay \emph{bounded} by making use of future information. However, the diversion polices of \cite{SSX12} require the length of the lookahead window to grow to infinity at a non-trivial rate in the heavy-traffic regime, and it remained open  whether substantial performance improvement could still be achieved with \emph{less} future information.

We resolve this question to a large extent by establishing an asymptotically tight lower bound on how much future information is {necessary} to achieve superior performance, which matches the upper bound of \cite{SSX12} up to a constant multiplicative factor. Our result hence demonstrates that the system's heavy-traffic delay performance is highly sensitive to the amount of future information available. Our proof is based on analyzing certain excursion probabilities of the input sample paths, and exploiting a connection between a policy's diversion decisions and subsequent server idling, which may be of independent interest for related dynamic resource allocation problems.}



\maketitle 


\section{Introduction}

Recently, there have been substantial interests in developing {forecasting systems} and {predictive models} across various application domains, which enable a system manager to obtain (partial) information of \emph{future} inputs, and thus allow for more efficient decision making or resource allocation. Examples of these systems include advanced ordering in supply chains (\cite{FR96}), appointment booking for elective surgeries (\cite{KH02}), and mechanisms for predicting future hospital visits (\cite{wardon_emj09,sun_09}). Because acquiring  accurate predictions can often involve additional infrastructural investments and operational complexities, it is a  natural question to ask  \emph{how useful} such predictive information can be, in terms of its ability in improving system performance beyond what can be achieved by the more conventional way of \emph{online} decision making, which does not take predictive information into account.

\begin{figure}[h]
\centering
\includegraphics[scale=1.1]{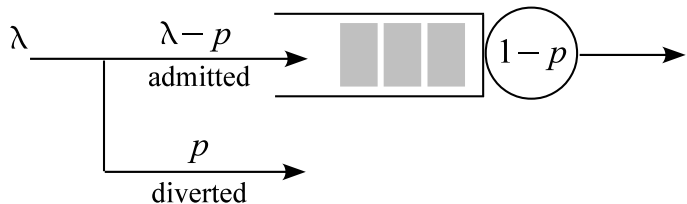}
\caption{An illustration of the queueing admission control problem.}
\label{fig:model}
\end{figure}

In a recent paper, \cite{SSX12} initiated an investigation along this direction in a class of {queueing admission control} problems, illustrated in Figure \ref{fig:model}.  An {overloaded} queue with service rate $1-p$ receives incoming jobs at rate $\lambda\in (1-p,1)$, and the system manager is allowed to \emph{divert} incoming jobs up to a rate of $p$, with the objective of minimizing the time-average queueing delay among the admitted jobs. The system manager has access to a \emph{lookahead window} of length $\wlam$, within which the realizations of future arrivals and service availability are revealed. The online version of the problem, with $\wlam=0$, is a classical queueing model that has been studied in various contexts related to congestion control (\cite{Yech71,Sti85}).

A main message of \cite{SSX12} is that one can drastically reduce queueing delay with  a \emph{sufficient} amount of future information. In particular, there exists $\ch>0$, such that if the length of the lookahead window satisfies
\begin{equation}
\wlam \geq \ch\ln \frac{1}{1-\lambda}, 
\label{eq:ssx12Req}
\end{equation}
then there exists a sequence of diversion policies, so that the resultant delay will stay \emph{bounded} in the heavy-traffic regime of $\lambda \to 1$. In sharp contrast, when no future information is available, the delay under an {optimal} online policy will diverge to \emph{infinity}, as $\lambda \to 1$. 

However, the requirement on the length of the lookahead window, as in Eq.~\eqref{eq:ssx12Req}, means that the superior delay performance achieved by  \cite{SSX12} comes  at the expense of a non-trivial amount of predictive power. Therefore, it remains to determine whether one could use much less future information and still achieve a significant performance improvement over an optimal online policy. This question is of practical importance, because  a larger amount of future information often requires more sophisticated  predictive models and computational infrastructures, which can be costly, if not impossible, to build and operate. 

The main contribution of the present paper is to provide a negative answer to above question, by showing that there exists a positive constant, $\cl$, such that if $\wlam$ scales \emph{slower} than $\cl\ln\frac{1}{1-\lambda}$ as $\lambda \to 1$, then the resulting delay performance can  be \emph{no better} than that of an optimal online policy by more than a constant factor. As a by-product of our result, an interesting ``conservation law'' is established, which suggests that {delay} and {future information} are, in some sense, ``exchangeable'' quantities (see discussions in Section \ref{sec:implication}). 

Despite having identical modeling assumptions, our proof techniques are quite different from those employed by \cite{SSX12}. The core of our arguments hinges upon a relationship between {diversions} and \emph{future idling} of the server, evaluated over certain subset of  input sample paths. This relationship is then used in conjunction with the excursion probabilities of a transition random walk to demonstrate that the system manager \emph{must} maintain a relatively large queue length, when the amount of future information is limited. We believe that this line of arguments is fairly robust to changes in modeling assumptions, and can be generalized, in other dynamic resource allocation problems, to proving lower bounds for the amount of information necessary in achieving desirable performance. 

\subsection{Organization} The remainder of the paper is organized as follows. In Section \ref{sec:res}, we state our main result, Theorem \ref{thm:lowerbound}, and contrast it with the prior results of \cite{SSX12}. In the same section, we discuss several implications of the theorem (Section \ref{sec:implication}), as well as connections of our work to the literature (Section \ref{sec:relWork}). Section \ref{sec:model} describes the modeling assumptions in more details, and introduces the necessary mathematical formalism. The proof of Theorem \ref{thm:lowerbound} is given in Section \ref{sec:proof}, with an outline of the proof ideas provided at the beginning of the section. We conclude the paper in Section \ref{sec:conclusion} and examine potential directions for future research.

\section{Main Result}
\label{sec:res}

\emph{Review of Prior Results}. We begin by informally reviewing the system model in \cite{SSX12}, which will be described in detail in Section \ref{sec:model}. The admission control problem runs in continuous time, and is characterized by three parameters: $\lambda$, $p$, and $W_\lambda$. An illustration of the system model is given in Figure \ref{fig:model}. 
\begin{enumerate}
\item  Jobs arrives to the system at the rate of $\lambda$, where $\lambda \in (0,1)$. There is a single server which processes jobs at the rate of $1-p$, where $p$ is a \emph{fixed} constant in $(0,1)$. It is assumed that the system is operating in the \emph{overload} regime, with $\lambda > 1-p$.

\item Upon each job's arrival, the system manager decides whether the job is to be admitted or diverted. If admitted, the job queues up in an (infinite) buffer until it is processed by the server, and if diverted, it leaves the system immediately. The goal of the system manager is to choose a \emph{diversion policy} that minimizes the time-average queue length induced by the admitted jobs, subject to the constraint that the infinite-horizon time-average rate of diversion does not exceed $p$. 

We will be primarily interested in the \emph{heavy-traffic regime} of $\lambda \to 1$,  where the post-diversion arrival rate approaches the server capacity of $1-p$, assuming that the system manager diverts at the maximum allowable rate of $p$. Note that by Little's Law, the time-average queue length is equal to the time-average queueing delay multiplied by the post-diversion arrival rate of $\lambda-p$. In the limit of $\lambda \to 1$, the two quantities will differ only by a multiplicative constant of $1-p$. Therefore, from this point on, we will focus on the time-average queue length as the performance metric, with the understanding that an analogous statement will hold for delay as well. 

\item The system manager has access to {information about the future}, which takes the form of a \emph{lookahead window} of length $\wlam$: at time $t$, the times of arrivals and service availability within the interval $[t,t+\wlam]$ are revealed to the system manager\footnote{Depending on the application, one can think of the lookahead window as being provided by some external oracle, or a predictive model that has access to side information.}. The case of $\wlam=0$ will be referred to the \emph{online} problem, since the system manager does not have access to any future information. 
\end{enumerate}

Denote by $\mathcal{Q}(\pi,\lambda,\wlam)$ the time-average queue length under the diversion policy $\pi$, given arrival rate $\lambda$ and a lookahead window of length $\wlam$. Let $\qop$ be the time-average queue length under an \emph{optimal} diversion policy (assuming such optimal policies exist), with
\begin{equation}
\qop = \min_{\pi} \mathcal{Q}(\pi,\lambda,\wlam), 
\end{equation}
It is shown in  \cite{SSX12} that a finite amount of lookahead into the future is sufficient to yield significant delay improvement over an online policy. In particular, fixing $p\in (0,1)$, they show that the optimal average queue length for an online policy diverges to {infinity} in the heavy-traffic regime, with
\begin{equation}
\mathcal{Q}^*(\lambda, 0) \sim \log_{\frac{1}{1-p}}\frac{1}{1-\lambda}, \quad \mbox{as $\lambda \to 1$}. 
\label{eq:onlineScale}
\end{equation}
In sharp contrast, there exists a positive constant $\ch$, whose value can depend on $p$, so that if 
\begin{equation}
 \wlam \geq \ch \ln \frac{1}{1-\lambda},
 \label{eq:futUpper}
\end{equation} 
for all $\lambda$ sufficiently close to $1$, then the optimal average queue length converges to a {finite} constant in the heavy-traffic regime: 
\begin{equation}
\qop \to \frac{1-p}{p}, \quad \mbox{as $\lambda \to 1$}. 
\label{eq:futureUpper}
\end{equation}

A main open question posed by \cite{SSX12} is whether significant performance gain over the online policy can still be achieved under much less future information. It is conjectured that if $\wlam = o\left(\ln \frac{1}{1-\lambda}\right)$, then the average queue length will necessarily diverge to infinity in the heavy-traffic limit (Conjecture $1$, \cite{SSX12}). In other words, a sufficient amount of future information may be \emph{essential} in achieving superior delay performance. 

\emph{Our Result}. The main result of this paper confirms, and strengthens, this conjecture of \cite{SSX12}. We show that if the amount of future information is \emph{insufficient} even by a constant factor, then not only will the delay be infinite in the heavy-traffic regime, but the delay scaling will essentially be \emph{no better} than that of an online policy. Specifically, we have the following theorem.

\begin{theorem}[Necessity of Future Information]
\label{thm:lowerbound}Fix $p\in (0,1)$. There exist $\cl> 0$ and $\tilde{\lambda} \in (1-p,1)$, so that if 
\begin{equation}
 \wlam \leq  \cl \ln\frac{1}{1-\lambda},\quad \forall \lambda \in (\tilde{\lambda},1),
 \end{equation}
then\footnote{The notation $f(x) = \Theta(g(x))$, as $x \to 1$,  represents the statement that, for any sequence $x_n\to 1$, we have $0<\liminf_{n \to \infty}{{f(x_n)}/{g(x_n)}}\leq \limsup_{n\to \infty}{{f(x_n)}/{g(x_n)}} < \infty$. }
\begin{equation}
\label{eq:qoplb}
 \qop = \Theta\left( \ln \frac{1}{1-\lambda}\right), \quad \mbox{as $\lambda \to 1$}.
\end{equation}
\end{theorem}

\begin{figure}[h]
\vspace{-8pt}
\begin{center}
\includegraphics[scale=1.15]{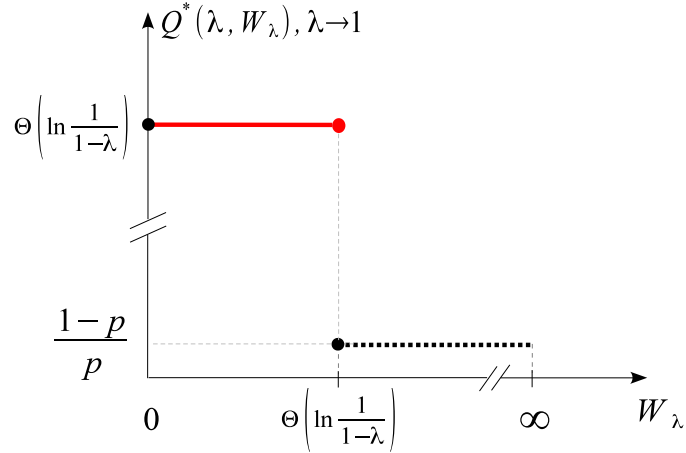}
\caption{Impact of future information on the effectiveness of admission control, in the heavy-traffic regime of $\lambda \to 1$. The solid red segment corresponds to the regime established by this paper, where $W_\lambda \lleq c_l\ln\frac{1}{1-\lambda} $ (Theorem \ref{thm:lowerbound}), and the dotted black segment corresponds to the regime established by \cite{SSX12}, where $W_\lambda \ggeq c_h\ln\frac{1}{1-\lambda}$ (Eq.~\eqref{eq:futUpper} and \eqref{eq:futureUpper} in the current paper). The case of $W_\lambda =0$ is covered by either paper. }
\label{fig:phaseTrans}
\end{center}
\end{figure}

Together with the results of \cite{SSX12}, Theorem \ref{thm:lowerbound} suggests that the performance of the admission control problem \emph{depends critically}  on the amount of future information available, and in particular, on how the length of the lookahead window, $\wlam$, scales relative to the watershed of $\Theta\left(\ln \frac{1}{1-\lambda}\right)$. A graphical illustration of Theorem \ref{thm:lowerbound}, with a comparison to the results of \cite{SSX12}, is provided in Figure \ref{fig:phaseTrans}. 

The proof of Theorem \ref{thm:lowerbound} is given in Section \ref{sec:proof}. It is worth noting that our proof techniques are quite different from those employed by \cite{SSX12}. In fact, they are somewhat ``{dual}'' to each other:  the earlier achievability result (Eq.~\eqref{eq:futureUpper})  was proved by analyzing the distribution of the lengths of {busy periods} associated with the queue length process (a property in \emph{time}), whereas the core of our arguments relies on the excursion properties of a transient random walk (a  property in \emph{space}).  

\subsection{Implications of Theorem \ref{thm:lowerbound}
}

There are several interesting implications of Theorem \ref{thm:lowerbound}. First, by virtue of being a lower bound for the case where the decision maker is given the exact realizations of  future input, Theorem \ref{thm:lowerbound} automatically extends to settings where predictions can be {noisy} or {corrupted}, as is typically the case in practical applications.   

Theorem \ref{thm:lowerbound} also implies an interesting ``conservation law'' between delay and future information: from Eqs.~\eqref{eq:onlineScale} through \eqref{eq:qoplb}, we see that the {sum} of  $\qop$ and $\wlam$ must be of order $\mathbf{\Omega}\left( \ln \frac{1}{1-\lambda}\right)$, as $\lambda \to 1$. In a rough sense, this is because the {same} type of stochastic discrepancies in the input processes, which  necessitate large queueing delays in the heavy-traffic when future information is limited, also determines how much lookahead is required in order to achieve a bounded delay. Even though such conservation seems to suggest that there is no ``free lunch'' to be had, the ability to understand and make such trade-offs can still be useful, because depending on the application, future information may be significantly less costly than delay, or \emph{vice versa}. 

From an operational point of view, although Theorem \ref{thm:lowerbound} invalidates the usefulness of future information in certain regimes, it is nevertheless reassuring to know that a simple online policy could do almost as well as any sophisticated prediction-guided policies, even when the amount of predictive information available grows as the traffic intensity increases. Moreover, the theorem does not rule out the possibility of having meaningful prediction-guided policies when future information is limited; it only implies that our search  in such scenarios should aim at more moderate, \emph{constant factor} performance improvements over online policies. In fact, numerical results in \cite{XuChan14} on a similar admission control model suggest that sizable performance gains can still be achievable, even with limited and noisy predictive information. 
\label{sec:implication}

\subsection{Related Work}
\label{sec:relWork}

In terms of modeling assumptions, our setup is identical to that of \cite{SSX12}, and hence we refer the reader to \cite{SSX12} for a review of the model's connections with the literature on classical Markov admission control problems and competitive analysis. 
The model is also related to a multi-server system with partial resource pooling (cf.~\cite{TX12}); the reader is referred to Chapter 7 of \cite{xu2014power} for more details. In addition, \cite{XuChan14} examines the model's relevance in the context of reducing waiting times at emergency departments. 

Our result can be viewed as a generalization of the {Markov} optimal admission control problem that has been studied in the literature (\cite{Sti85}), and it is interesting to contrast some of the differences in analytical approaches. Optimal policies in the Markov setting ($\wlam=0$) are  known to often admit a \emph{threshold} (or control-limit) form, where a diversion is made only if the current queue length reaches a fixed threshold. To prove the optimality of these policies, one would typically analyze the Bellman equations of the corresponding Markov decision process (MDP) in order to establish a set of {monotonicity} properties in the policy space, e.g., that the cost-to-go function for a threshold policy would be dominated by policies that divert with non-zero probabilities when the queue is small (c.f.~\cite{Yech71}). Successive applications of such monotonicity properties will then narrow the policy space down to only those with a threshold form. 

Unfortunately, these arguments employed in the Markov setting do not seem to carry over easily when the lookahead window is taken into account. While our setting can still be cast as an MDP by incorporating the lookahead window into the state space, the structure of the state space is now considerably more complex (and increasingly so, as $\wlam  \to \infty$), and it is not so clear as to whether any monotonicity property continues to hold. Our proof techniques circumvent this complexity by focusing on the ``macroscopic'' sample-path characteristics of the system, instead of the more refined details of the Bellman equations. As a trade-off, our analysis is more ``coarse'' by nature, and it provides neither a characterization of the multiplicative \emph{constant} in the delay scaling, nor a concrete diversion policy that achieves the lower bound of the necessary amount of future information (which, fortunately, has already been given in \cite{SSX12}). 

Our work is also similar in spirit to the techniques of information relaxation and path-wise optimization for MDPs (\cite{Rogers07,Brown10,Farias12}). In this case, one considers an relaxed version of the original MDP, where the decision maker has access to realizations of the future input sample paths. This relaxed problem is often  simpler to solve and simulate than the original stochastic optimization problem, and hence can be used, for instance, as a performance benchmark for evaluating heuristic policies. Our work is different from this literature in several aspects. Most notably, we focus on rigorously understanding the stochastic dynamics involved in the relaxed problem with future information, and how performance scales with respect to the length of the lookahead window, as opposed to using the relaxed problem to approximate the performance of an optimal online policy, which is well understood in our setting.

\section{Model and Notation}
\label{sec:model}

We now present the mathematical formalism and modeling assumptions that will be used throughout the remainder of the paper. An illustration of the system is given in Figure \ref{fig:model}. 

\emph{System Dynamics}. The system runs in continuous time, indexed by $t\in \rp$. There is a \emph{queue} with infinite waiting room, whose length at time $t$ is denoted by $Q(t)$.  The input to the system consists of two independent Poisson processes: 
\begin{enumerate}
\item $\arr$, with rate $\lambda$, which corresponds to the \emph{arrival} of jobs; 
\item $\ser$, with rate $1-p$, which corresponds to the generation of \emph{service tokens}.
\end{enumerate}
When an \emph{event} occurs in $\arr$ at time $t$, we say that a job has arrived to the system, and the value of $Q(t)$ is incremented by $1$, if the job is ``admitted'' (see below for the description of admission policies). Similarly, when an event occurs in the process $\ser$ at time $t$, we say that a service token is generated, and the value of $Q(t)$ is decremented by $1$, if $Q(t)>0$, and remains at $0$, otherwise.\footnote{The generation of a service token at time $t$ can be thought of as the server being able to fetch a new job from the queue at time $t$. As such, the service token model attributes the randomness in processing times to an external source, which does not depend on the identities of the jobs. It can be shown that, in the online setting, the service token model is equivalent to the more conventional assumption of exponentially distributed job sizes, though such equivalence is generally not true when future information is taken into account. The reader is referred to Page 9 of \cite{SSX12}, and the references therein, for more details on the service token model. }

For our purposes, it is more convenient to work with the sequence $\{(Z_n,R_n): n \in \N\}$, where 
\begin{equation}
Z_n = \mbox{time of the $n$th event in $\arr\cup \ser$},
\end{equation}
and $R_n$ encodes the type of the $n$th event, with
\begin{equation}
R_n = \left\{
\begin{array}{ll}
1,& \quad \mbox{if the $n$th event is in $\arr$ (arrival)}, \\
-1,& \quad  \mbox{if the $n$th event is in $\ser$ (service token)}.
\end{array}
\right.
\end{equation}
We will let $\{\cnt(t): t\in \rp\}$ be the counting process associated with $\{Z_n\}$, with 
\begin{equation}
	\cnt(t) =  \sup \{n \in \zp: Z_n \leq t\}, 
	\label{eq:cntDef}
\end{equation}
and denote by $S(s,t)$ the \emph{difference} between the numbers of arrival and services tokens in the interval $(s,t]$,
\begin{equation}
	S(s,t)= \sum_{\cnt(s)+1 \leq n \leq \cnt(t)} R_n.
	\label{eq:Sdef}
\end{equation}
Note that when $\lambda \neq 1-p$ the process $\{S(0,t): t\in \rp\}$ is a transient random walk, with 
\begin{equation}
\E(S(0,t)) = [\lambda-(1-p)]t. 
\end{equation}

\emph{Future Information}. The notion of future information  is captured by a \emph{lookahead window}. At any time $t$, the system manager has access to the {realization} of all events in $\arr\cup \ser$ in the interval $[t,t+\wlam]$. Throughout, we will denote by $\wlam$ the length of the lookahead window, under arrival rate $\lambda$. 

\emph{Admission Policies}. Upon arrival, each job is either \emph{admitted}, in which case it joins the queue, or \emph{diverted}, in which case it disappears from the system immediately. The role of a diversion policy, $\pi$, is to output a sequence of {diversion decisions} for all events, represented by the sequence of indicator variables, $\{H(n) : n\in \N\}$, where
\begin{equation}
	H(n) = \mathbb{I}\left\{\mbox{$R_n=1$, and $\pi$ chooses to divert at time $Z_n$}\right\}. 
\end{equation}
Given the form of future information, we will require  that the diversion policy be $(t+\wlam)$-causal, so that the decision made at time $t$ does not depend on any event after time $t+\wlam$.  A diversion policy is said to be  \emph{feasible}, if the resulting time-average rate of diversion is at most $p$, i.e.,
\begin{equation}
\label{eq:avgCost}
	\limsup_{N\to \infty}\frac{\lambda+1-p}{N} \E\left(\sum_{n=1}^{N}H(n)\right) \leq p. 
\end{equation}
where the constant $\lambda+1-p$ corresponds to the total rate of events in $\arr\cup \ser$. The objective of the decision maker is to choose a feasible policy, $\pi$, so as to \emph{minimize} the time-average queue length, defined by\footnote{Throughout, $f(x-)$  represents the limit $\lim_{y \uparrow x}f(y)$.}
\begin{equation}
\label{eq:avgQueue}
	\mathcal{Q}(\pi,\lambda,\wlam)= \limsup_{N\to \infty} \E\left(\frac{1}{N} \sum_{n=1}^N Q(Z_n-)\right).
\end{equation}

\subsection{Notation}
We will assume that all asymptotic expressions with respect to $\lambda$ are taken in the limit of $\lambda \to 1$.  We will use $f\ll g$ and  $f \lleq g$ to denote $f=o(g)$ and $f=\mathcal{O}(g)$, respectively. We will write $f \asleq g$ to mean that $f(x)\leq g(x)$ for all $x$ sufficiently closely to $1$, i.e., that there exists $y\in (0,1)$, such that $f(x)\leq g(x)$, for all $x \in (y,1)$. The expressions $f \gg$, $\ggeq$ and  $\asgeq g$ are defined analogously to their respective counterparts. When a statement is made concerning the limit ``as $x\to 1$'', without specifying the exact sequence with respect to which the limit is taken, it is understood that the statement should hold for any sequence, $\{x_n\}$, with $\lim_{n\to \infty}x_n=1$. The notation $X\stackrel{d}{=} Y$ means that the random variables $X$ and $Y$ have the same distribution.

\section{Proof of Theorem \ref{thm:lowerbound}} 
\label{sec:proof}

The remainder of the paper is devoted to the proof of Theorem \ref{thm:lowerbound}. We begin with a high-level summary of the main steps involved. First, we argue that there exists a stationary optimal policy,  which makes decisions only based on the current queue length and the content of the lookahead window. Furthermore, the queue length process under this stationary policy admits a well-defined steady-state distribution  (Section \ref{sec:stationary}). This stationarity will allow us to simplify the analysis by focusing  on the policy's actions over a finite time horizon.

We will prove Theorem \ref{thm:lowerbound} by contradiction, where we start by assuming that a small average queue length is indeed achievable under an optimal stationary policy, even with a small lookahead window, and later refute this assumption. Our main arguments are based on the identification of a set of \emph{base sample paths} (Section \ref{sec:baseSamp}), with the property that \emph{any} feasible policy must perform poorly over these sample paths, should the length of the lookahead window be too small. The stationarity property described earlier will then allow us to extend this argument to showing the policy's failure over the infinite time horizon. It is worth noting that the base sample paths are not ``typical,'' in the sense that their occurrences possess only vanishingly small probability, as $\lambda \to 1$.   This is because the failures of a policy under a small lookahead window are not caused by the average behavior of the inputs, but rather by some rare excursions of the  random walk $S(0,\cdot)$. Though occurring with small probabilities,  these excursions are in some sense unforeseeable under a small lookahead window, and their existence forces an optimal policy to be overly restrained in diverting jobs and hence yield a large average queue length.  

 To carry out the arguments using the base sample paths, we will exploit a key relationship between \emph{diversions} and \emph{server idling}. In particular, we will demonstrate that, without sufficient lookahead, if a {constant fraction} of the arrivals are diverted during a specific portion of a base sample path, it will inevitably result in excessive idling of the server not far away in the future, even as $\lambda \to 1$. However, such server idling cannot occur in the heavy-traffic limit, since the server must be fully utilized in order to ensure system stability. This reasoning then implies that any  policy that makes such diversions must be {infeasible}, or conversely, that any feasible policy must divert very few arrivals over these segments of the base sample paths (Proposition \ref{prop:nodel}). However, such conservatism comes at a cost, in that it leads to long episodes during which the queue length stays at a  high level (Proposition \ref{prop:D1lowb}).  We then argue that the frequent appearances of such ``bad''  episodes will result in a large average queue length in steady-state, which contradicts with our initial assumption and hence completes the proof of Theorem \ref{thm:lowerbound}.

\subsection{Preliminaries}
Without loss of generality, we will consider only the cases where the length of the lookahead window, $\wlam$, 
diverges to infinity in the heavy-traffic regime, i.e., 
\begin{equation}
\wlam \to \infty, \quad \mbox{as $\lambda \to 1$.}
\label{eq:Wlamassump}
\end{equation} 
To see why this is justified, note that because we can always achieve the same average queue length with a longer lookahead window,  the optimal average queue length $\qop$  must be monototically non-increasing in  $\wlam$. Therefore, any lower bound we obtain  on $\qop$ under the assumption of Eq.~\eqref{eq:Wlamassump} also applies to the case where $\wlam = \mathcal{O}(1)$. For simplicity of notation, we will drop the dependency on $W_\lambda$, and denote by $q_\lambda$ the optimal average queue length, 
\begin{equation}
\qlam=\qop, \quad \forall \lambda \in (1-p,1).
\end{equation}

\emph{Main Assumption}. We  will assume the validity of the following hypothesis throughout the remainder of the proof, which states that it is indeed possible to achieve a small delay as long as $\wlam$ is of order $\mathbf{\Omega}\left(\ln \frac{1}{1-\lambda}\right)$. As will be shown in Section \ref{sec:finishingProof}, invalidating this hypothesis will imply the lower bound in Theorem \ref{thm:lowerbound}. 

\begin{hypothesis}
\label{ass:main}
Fix $p\in (0,1)$. Suppose that $\wlam \ggeq \ln\frac{1}{1-\lambda}$, as $\lambda \to 1$. Then
\begin{equation}
\label{eq:contra1}
q_{\lambda} \ll \ln\frac{1}{1-\lambda}, \quad \mbox{as $n \to \infty$}. 
\end{equation}
\end{hypothesis}
Assuming the validity of Hypothesis \ref{ass:main}, it also follows that if $\wlam\ggeq \ln \frac{1}{1-\lambda}$, as $\lambda \to 1$, then
\begin{equation}
\qlam \ll \wlam, \quad \mbox{as $\lambda \to 1$.}
\label{eq:ass1a}
\end{equation}

\subsubsection{State Representation and Stationary Policies}
\label{sec:stationary}

We show in this section that there always exists an \emph{stationary} optimal policy that depends only on the {state}, which consists of the current queue length and content of the lookahead window. 

Since all diversion decisions are associated with events in $\arr\cup \ser$, it suffices to specify the nature of future information for the event times, $\{Z_n: n\in \N\}$. 
At $t=Z_n$, the \emph{content} of the lookahead window is defined to be the vector $F(n)=(F_k(n)\,:\, k \in \zp)$, where
\begin{equation}
	F_k(n) = (Z_{n+k}-Z_n \,,\, R_{n+k}), \quad  0\leq k\leq \cnt(Z_n+\wlam)-\cnt(Z_n).
	\label{eq:defFn}
\end{equation}
In other words, $F_k(n)$ specifies the time of the $k$th future event starting from the current time, $Z_n$, along with its type for all events within the lookahead window of length $\wlam$. For future events beyond the lookahead window which we have no access to, we simply set the value of $F_k(n)$ to zero: 
\begin{equation}
	F_k(n) = (0 \,,\, 0), \quad  k > \cnt(Z_n+\wlam)-\cnt(Z_n). 
\end{equation}

Recall that $Q(t)$ is the queue length at time $t$. Consider the sequence $\{X(n): n\in \N\}$, where
\begin{equation}
	X(n) = \left(Q(Z_n-), F(n)\right).
\end{equation}
From this point on, we will refer to $\left\{X(n): n\in \N\right\}$  as the \emph{states} of our system. 

\emph{Stationary Policies}. A diversion policy $\pi$ is \emph{stationary}, if its diversion decision at time $Z_n$ depends only on the state, $X(n)$, or formally, that
\begin{equation}
	\pb\left(H(n)=1 \,\big|\, X(n)\right) = \pb\left(H(n)=1 \,\Big|\, \left\{(Z_k,R_k)\right\}_{k=1}^{\cnt(Z_n+\wlam)}\right), \quad\mbox{a.s.}
\end{equation} 
A stationary policy, $\pi$, is \emph{stable}, if the evolution of $\{X(n): n\in \N\}$ under $\pi $ admits a well-defined steady-state distribution, $\gamma$, so that steady-state queue length and probability of diversion coincide with the time-average queue length and diversion rate, respectively, given that the initial condition, $X(1)$, is distributed according to $\gamma$. 

In our admission control problem, because the arrivals and service tokens are generated according to Poisson processes, future evolution of the system starting from $t=Z_n$ is independent conditional on the current state $X_n$ and diversion decision. As such, our problem can be cast as a discrete-time Markov decision process (MDP), with states $\{X_n: n\in \N\}$ and actions that correspond to the probabilities of diversion. Using existing results in the literature (c.f.~\cite{HGL03,GV11}), it can be shown that, for MDPs of this kind, there exists an optimal policy that is also stationary and stable. This is summarized in the following lemma, whose proof is given in Appendix \ref{app:lem:stationaryStable}. 
\begin{lemma}
\label{lem:stationaryStabel}
Fix any $p>0$, $\lambda\in (1-p, 1)$, and $\wlam>0$. The admission control problem admits a \emph{stable stationary} optimal policy, $\pi$, which achieves the minimum time-average queue length among all feasible diversion policies. 
\end{lemma}

In light of Lemma \ref{lem:stationaryStabel}, we will, in the remainder of the proof of Theorem \ref{thm:lowerbound}, focus on the family of stable stationary policies, which we will refer to simply as \emph{stationary policies}. Given a stationary policy, $\pi$, the resultant state sequence $\{X(n): n\in \N\}$ is a stationary Markov chain. Since we are interested in deriving a performance lower bound, we may assume that, at time $t=0$,  both the queue length and the content of the lookahead window are initialized according to the steady-state distributions, $\gamma$. In particular, we have that 
\begin{equation}
	\E\left(Q(t)\right)=\E\left(Q(0)\right)=\mathcal{Q}(\pi,\lambda, \wlam), \quad t\in \rp. 
\end{equation}
and, that 
\begin{equation}
\E(H(n)) = \E(H(1))=\limsup_{N\to \infty}\frac{\E\left(\sum_{n=1}^N H(n)\right)}{N}, \quad \forall n \in \N. 
\label{eq:Hstation}
\end{equation}

Define the process $\left\{L(t): t\in \rp\right\}$, where
\begin{equation}
	L(t) = \mathbb{I}\left\{Q(t)\leq 2\qlam\right\}, \quad t\in \rp.
\end{equation}
The following lemma will be useful. 
\begin{lemma} Fix $p\in (0,1)$. For all $\lambda\in (1-p,1)$, we have that
\begin{equation}
	\label{eq:visitFreqBound}
	\E (L(t)) = \pb\left(Q(0)\leq 2\qlam \right) \geq \frac{1}{2},\quad \forall t\in \rp, 
\end{equation}
under any optimal stationary policy. 
\end{lemma}
\bproof
The result follows from the stationarity of $Q(\cdot)$ and the Markov's inequality: 
\begin{equation}
	\E (L(t)) = \pb\left(Q(t)\leq 2\qlam \right)  = \pb\left(Q(0)\leq 2\qlam \right) = \pb\left(Q(0)\leq 2\E(Q(0)) \right) \geq \frac{1}{2}. 
\end{equation}
\qed
%

In the remainder of the proof, we will show that there exists $\cl>0$ such that if $\wlam \asleq \cl \ln\frac{1}{1-\lambda}$, then Eq.~\eqref{eq:visitFreqBound} cannot be true under any sequence of optimal stationary policies, unless $Q^*(\lambda, \wlam) \ggeq \ln\frac{1}{1-\lambda}$. This would invalidate Hypothesis \ref{ass:main}, which would in turn prove the lower bound on $\qop$ in Theorem \ref{thm:lowerbound}. 


\subsection{Base Sample Paths} 

We now describe the construction of a set of base sample paths which will serve as the basis of our subsequent analysis. In later sections, we will show that, roughly speaking, the non-negligible chance of occurrence of such sample paths will ``force'' any feasible policy to be overly conservative in diverting jobs, should $\wlam$ be too small.

\label{sec:baseSamp}
\begin{figure}[h]
\begin{center}
\includegraphics[scale=.65]{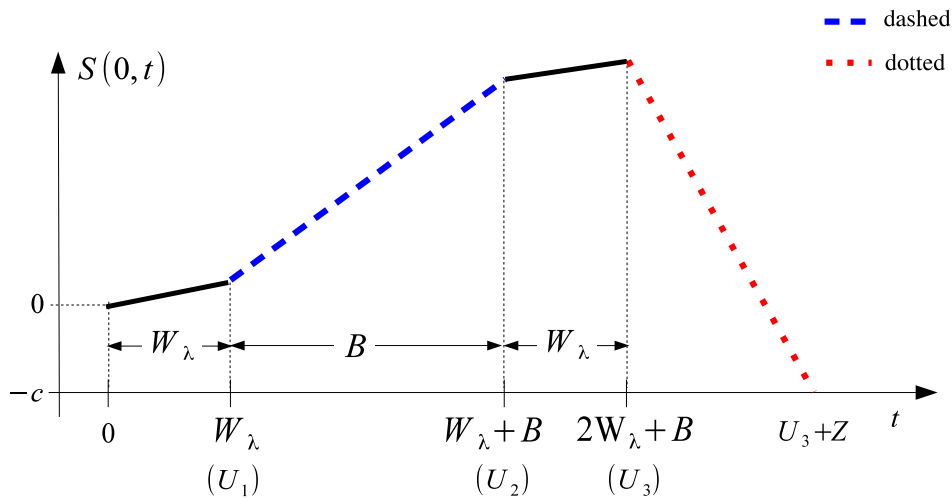}
\caption{This figure illustrates the ``macroscopic'' behavior of the base sample paths.  The dashed blue segment between $\wlam$ and $\wlam+B$ represents a period of sustained upward drift of $S(0,\cdot)$, and the dotted red segment starting at $2\wlam+B$ represents a downward drift.  The two solid black segments, each with length equal to that of the lookahead window, serve as a ``buffer'', ensuring that the actions of the diversion policy before the segment are {independent} from the evolution of $S(0,\cdot)$ afterwards. }
\label{fig:samplePath}
\end{center}
\end{figure}

Let $B\in \rp$ be a quantity whose value will be specified in the sequel. We define the following time markers, whose positions relative to each other are illustrated in Figure \ref{fig:samplePath}.  
\begin{align}
	U_1 &= \wlam, \nnb\\
	U_2 &= U_1+B=\wlam+B, \nnb \\
	U_3 &= U_2+\wlam=2\wlam+B.\nnb
\end{align}
The set of base sample paths is defined as the intersection of the events $\eva$ through $\eve$, described as follows. Let $\epsilon, \zeta$ and $\phi$ be positive constants. 
\begin{enumerate}
\item Event $\eva$, parameterized by $\epsilon$ and $\zeta$, says that the sample path of $S(0,\cdot)$ stays close to its expected behavior during the interval $(U_1, U_2]$: 
\begin{equation}
\eva = \left\{\left| S\left(U_1,t\right) - [\lambda-(1-p)] t\right| \leq \epsilon t + \zeta , \mbox{ for all } t\in (U_1, U_2] \right\},
\label{eq:E2def}
\end{equation}
When $\epsilon$ is small, this implies that $S(0,\cdot)$ undergoes a consistent \emph{upward} drift during $(U_1,U_2]$. Event $\eva$ is illustrated by the dashed blue line segment in Figure \ref{fig:samplePath}.
 
\item Event $\evb$ says that the queue length at $t=0$ is not too large compared to the optimal average queue length,
\begin{equation}
\evb = \left\{Q(0)\leq 6\qlam\right\}. 
\label{eq:E3def}
\end{equation}

\item The events $\evc$ and $\evd$ put some restriction on the amount of upward excursion of $S(0, \cdot)$ during the intervals $(0,U_1]$ and $(U_2, U_3]$, respectively, 
\begin{equation}
\evc = \left\{S(0,U_1) \leq 2W_\lambda \right\} ,
\label{eq:E4def}
\end{equation}
\begin{equation}
\evd = \left\{S(U_2,U_3) \leq 2W_\lambda \right\},
\label{eq:E5def}
\end{equation}
The main purpose of $\evc$ and $\evd$ is to serve as ``buffers'' to induce certain independence property, which will be useful for subsequent analysis: since the lengths of $(0,U_1]$ and $(U_2, U_3]$ are both equal to that of the lookahead window, the actions of the diversion policy before each interval are \emph{independent} from the evolution of $S(0,\cdot)$ after it. The two events are illustrated by the black line segments in Figure \ref{fig:samplePath}.

\item Finally, the event $\eve$ says that $S(0,\cdot)$ will undergo a substantial \emph{downward} excursion soon after $U_3$, as is illustrated by the dotted red line segment in Figure \ref{fig:samplePath}. Let $Z$ be the stopping time 
\begin{equation}
\label{eq:Zdef}
	Z = \inf\left\{z \in \rp: S\left(U_3, U_3+z \right) < -\left[6\qlam+ [\lambda-(1-p)-\epsilon]B+\zeta+ 4\wlam \right] \right\},
\end{equation}
and $\eve$ is defined by putting an upper bound on $Z$: 
\begin{equation}
\eve =\left \{ Z \leq \phi\wlam \right\}.
\label{eq:E6def}
\end{equation}
The right-hand-side of the inequality in the definition of $Z$ was chosen so that, conditional on the joint occurrence of $\eva$ through $\evd$, a downward excursion in $S(0,\cdot)$ of such magnitude is guaranteed to \emph{deplete} the queue by time $U_3+Z$. As will become clearer in the next section, this depletion will help us connect diversions to future idling of the server. 
\end{enumerate}
Note that the events $\eva$, $\evc$, $\evd$ and $\eve$ concern the input sample path $S(0,\cdot)$ only, and are independent of the diversion policies, while $\evb$ also depends on the choice of diversion policy. 

Having described the events that together characterize the base sample paths, we next illustrate some of their statistical properties. The first lemma shows that the events $\eva$ through $\evd$ can occur with fairly high probabilities. The proof is given in Appendix \ref{app:lem:eveprop}. 

\begin{lemma} 
\label{lem:eveprop}
\begin{enumerate}
\item Fix $\epsilon >0$. For all $\theta \in (0,1)$, there exists $\zeta>0$, so that for all $\lambda > 1-\frac{1}{2}p$, 
\begin{equation}
\label{eq:E2isbig}
\inf_{B\geq 0} \pb\left(\eva\right) \geq \theta. 
\end{equation}
\item  Under optimal stationary policies, $\pb\left(\evb\right) = \pb(Q(0)\leq 6\qlam)\geq \frac{5}{6}$, for all $\lambda\in (1-p,1)$. 
\item $\lim_{\lambda \to 1}\pb\left(\evc\right)=\lim_{\lambda \to 1}\pb\left(\evd\right)=1$.
\end{enumerate}
\end{lemma}

The next lemma shows that the event $\eve$ occurs with a small yet non-negligible probability. The proof is given in Appendix \ref{app:lem:PoissonLDP}. 

\begin{lemma}
\label{lem:Eve}
Fix $k, \phi, \zeta>0$, and $\epsilon \in \left(0\, ,\, \min\{\zeta, \lambda-(1-p)\}\right)$. Suppose that $B=k\wlam$, and $\qlam \ll \wlam$, as $\lambda\to 1$.  There exists $\gamma>0$, such that
\begin{equation}
	\pb\left(\eve\right) \ggeq \exp\left(-\gamma\wlam \right), \quad \mbox{as $\lambda \to 1$}. 
\end{equation} 
\end{lemma}

Finally, the following independence properties among the events will be useful. The proof is given in Appendix \ref{app:lem:eveIndp}. 

\begin{lemma} 
\label{lem:eveIndp} Fixing a feasible diversion policy, the following holds. 
\begin{enumerate} 
\item The events $\eva, \evc, \evd$ and $\eve$ are mutually independent. 
\item The event $\evb$ is independent of $\eva, \evd$ and $\eve$, but not necessarily $\evc$. 
\item Denote by $Y$ the number of diversions in the interval  $\left(U_1,U_2\right]$, i.e.,
\begin{equation}
Y = \sum_{\cnt(U_1)+1 \leq n\leq \cnt(U_2)} H(n). 
\label{eq:Ydef}
\end{equation}
Then $Y$ is independent of $\eve$. 
\end{enumerate}
\end{lemma}

\subsection{From Diversions to Server Idling}

The goal of this subsection is to show that, if $\wlam$ is small, then the number of diversions made during the the interval $(U_1,U_2]$, i.e., the random variable $Y$ (Eq.~\eqref{eq:Ydef}), must also be appropriately small, under any optimal stationary policy. To achieve this, we will exploit a connection between $Y$ and the idling of the server at a later time. 

The intuition is perhaps best seen pictorially, as depicted in Figure \ref{fig:samplePath}. Conditional on the occurrence of the events $\eva$ through $\eve$, and suppose {no} diversion has been made, the queue length process $Q(t)$ would have ``followed'' the trajectory depicted in the figure and reached {zero} by time $U_3+\phi \wlam$. Suppose now that a {large} number of diversions are made during the interval $(U_1, U_2]$ (dashed line segment in blue), the depletion of the queue implies that there must be an extended period of server idling prior to $U_3+\phi\wlam$. Such idling, if it persists even as $\lambda \to 1$, can be problematic and will be shown to contradict the feasibility of the diversion policy. This in turn implies that the number of diversions in  $(U_1, U_2]$ must be small. 

The next proposition is the main result of this subsection, which formalizes the above intuition. There is, however, one adjustment: as opposed to conditioning on all five events, which has vanishingly small probability due to the presence of $\eve$, we will condition only on $\eva$ and $\evb$, which occur with high probability. To do so, we will exploit several independence properties among the events, as in Lemma \ref{lem:eveIndp}, and show that the impact of $S(0,\cdot)$'s downward excursion described by $\eve$ is unavoidable when $\wlam$ is too small, even without explicitly conditioning on $\eve$. 

\begin{proposition}
\label{prop:nodel}
Fix $k>0$, and let $B=k\wlam$. There exists $c>0$, so that if 
\begin{equation}
	\wlam \asleq c\ln\frac{1}{1-\lambda}, \quad \mbox{as $\lambda \to 1$},
\end{equation}
then for every $\tau >0$,
\begin{equation}
\label{eq:Yis0}
	\lim_{\lambda \to 1} \pb\left(Y\geq \tau B \,\big|\, \eva \cap \evb \right) = 0,
\end{equation}
under any sequence of optimal stationary policies, where $Y$ is the number of diversions during  $(U_1, U_2]$, defined in Eq.~\eqref{eq:Ydef}. 
\end{proposition}

\bproof
We say that a service token generated at time $t$ is \emph{wasted}, if there is currently no job in the queue, i.e., $Q(t)=0$. Let $\{\jc(t): t\in \rp\}$ be the counting process of wasted service tokens, where
\begin{equation}
\jc(t) = \mbox{\# of wasted service tokens in $[0,t]$}. 
\end{equation}

For the sake of contradiction, assume the following is true: if  $\wlam \ggeq \ln\frac{1}{1-\lambda}$ as $\lambda \to 1$, then there exist $\tau>0$, and a sequence of optimal stationary policies, $\{\pi_\lambda\}$, under which
\begin{equation}
\label{eq:Ycontra1}
\liminf_{\lambda \to 1}\pb\left(Y\geq \tau B \,\big|\, \eva \cap \evb  \right)=q>0. 
\end{equation}
The following lemma is a key ingredient to the proof, which says that the number of wasted tokens must be substantial. The proof is based on the intuition explained in the passages above Proposition \ref{prop:nodel}, and is given in Appendix \ref{app:lem:EJLB}. 
\begin{lemma} 
\label{lem:EJLB}
Fix $k>0$, and let $B=k\wlam$. Suppose Eq.~\eqref{eq:Ycontra1} is true for some sequence of optimal stationary policies, $\{\pi_\lambda\}$. Then  there exist  $a,\gamma>0$ (whose values can depend on $k$) such that
\begin{equation}
\label{eq:lb_waste1}
	 \E\left(\jc(a\wlam)\right) \ggeq \wlam\exp\left(-\gamma \wlam\right) , 
\end{equation}
as $\lambda \to 1$, under $\{\pi_\lambda\}$. 
\end{lemma}

Consider an optimal stationary policy. Denote by $\hc(t)$ the counting process representing the number of diversions in $[0,t]$, i.e.,
\begin{equation}
\hc(t) = \sum_{n=1}^{\cnt(t)} H(n).
\end{equation}
By the stationarity of $\{H(n): n \in \N\}$ (Eq.~\eqref{eq:Hstation}) and definition of $\cnt(t)$ (Eq.~\eqref{eq:cntDef}), it is not difficult to show that, for all $t > 0$, 
\begin{align}
\frac{\E(\hc(t))}{t} =& \frac{1}{t}\E\left(\sum_{n=1}^{\cnt(t)} H(n)\right)= (\lambda+1-p)\E(H(1)) \nnb\\
=& \limsup_{N\to \infty}\frac{(\lambda+1-p)\E\left(\sum_{n=1}^N H(n)\right)}{N}. 
\label{eq:HExpContTime}
\end{align}
By definition, we have that
\begin{equation}
Q(t) = Q(0)+S(0,t)+\jc(t)-\hc(t), \quad \forall t >0. 
\end{equation}
Taking expectation on both sides of the above equation, and letting $t=a\wlam$, where $a$ is given as in Lemma \ref{lem:EJLB}, we have that
\begin{align}
&\frac{\E(\hc(a\wlam)) }{a\wlam}-p \nnb\\
=& \frac{1}{a\wlam}\left(\E\left(S\left(0,a\wlam\right)\right)+\E\left(\jc(a\wlam)\right) + \E(Q(0))-\E(Q(a\wlam))\right)-p \nnb\\
\sk{a}{=} & [\lambda-(1-p)]-p + \frac{1}{a\wlam}\E\left(\jc(a\wlam)\right) \nnb\\
\sk{b}{\ggeq}  &   (\lambda-1) + \frac{1}{a\wlam}\wlam\exp\left(-\gamma\wlam\right) \nnb\\
{\ggeq}& \exp\left(-\gamma\wlam\right) - (1-\lambda),
\label{eq:Havg1}
\end{align}
where $\gamma$ is given in Lemma \ref{lem:EJLB}. Step $(a)$ follows from the fact that $\E(Q(0))=\E(Q(a\wlam))$ by the stationarity of $Q(\cdot)$, and $(b)$ from Eq.~\eqref{eq:lb_waste1}.

Letting $\wlam = c \ln\frac{1}{1-\lambda}$, with   $c = {1}/{2\gamma}$, we have that 
\begin{equation}
\exp\left(-\gamma\wlam\right) \ggeq \sqrt{1-\lambda},  \quad \mbox{as $\lambda \to 1$}. 
\label{eq:Havg2}
\end{equation}
Combining Eqs.~\eqref{eq:Havg1} and \eqref{eq:Havg2}, we have that
\begin{equation}
\frac{\E(\hc(a\wlam)) }{a\wlam} - p \ggeq \sqrt{1-\lambda}-(1-\lambda) \ggeq  \sqrt{1-\lambda}, \quad \mbox{as $\lambda \to 1$}. 
\end{equation}
In particular, this implies that there exists $\lambda'\in (1-p,1)$, such that 
\begin{equation}
\frac{\E(\hc(a\wlam)) }{a\wlam}  > p, \quad \forall \lambda \in (\lambda',1).
\label{eq:hcContra}
\end{equation}
Since the stationary diversion policies we consider are feasible, we must have that 
\begin{equation}
\frac{\E(\hc(t))}{t} \stackrel{(a)}{=}\limsup_{N\to \infty}\frac{(\lambda+1-p)\E\left(\sum_{n=1}^N H(n)\right)}{N}  \stackrel{(b)}{\leq} p, 
\end{equation}
for all $\lambda \in (1-p,1)$, and $t>0$, where  $(a)$ and $(b)$ follow from Eqs.~\eqref{eq:HExpContTime} and \eqref{eq:avgCost}, respectively.  This leads to a contradiction with Eq.~\eqref{eq:hcContra}, which invalidates the assumption made in Eq.~\eqref{eq:Ycontra1}, and hence proves Proposition \ref{prop:nodel}.  $\qed$

\subsection{Consequences of Too Few Diversions}
Proposition \ref{prop:nodel} tells us that, under optimal stationary policies, the number of diversions in $(U_1, U_2]$ must be small when $\wlam$ is small. Building on this observation, we now focus on policies that divert ``very few'' jobs during $(U_1, U_2]$,  i.e., with  $Y$ scaling sub-linearly with respect to $B$, and show that they will necessarily lead to a large expected queue length in steady-state. The following proposition is the main result of this subsection. 

\begin{proposition}
\label{prop:D1lowb}
Fix $p\in (0,1)$. There exists $\cl> 0$, so that if 
\begin{equation}
	\wlam  \asleq \cl \ln\frac{1}{1-\lambda}, \quad \mbox{as $\lambda \to 1$},
\end{equation}
then
\begin{equation} 
\label{eq:ED1lim}
\limsup_{\lambda \to 1} \E\left(L(0)\right) \leq \frac{1}{3}. 
\end{equation}
under any sequence of optimal stationary policies. 
\end{proposition}


\bproof 
We will assume that $B=k\wlam$, with $k=24$, and that $\wlam \asleq c_l\ln \frac{1}{1-\lambda}$, where $\cl$ is equal to the constant $c$ in Proposition \ref{prop:nodel} for the corresponding value of $k$. 

Consider an optimal stationary policy, with a resultant average queue length of $\qlam$. We will prove the claim by showing that if $\eva\cap \evb$ occurs \emph{and} the number of diversions made in $(U_1,U_2]$ is small (cf.~Eq.~\eqref{eq:Yis0}), then, for a ``long time'' after $U_1$, the queue length will stay at a high level (i.e., $Q(t)>2\qlam$). Recall that $Y$ is the number of diversions made during the period $(U_1, U_2]$. We have the following inequality, derived from the queueing dynamics: 
\begin{equation}
\label{eq:QtU1U2}
Q(t) \geq Q(U_1) + S(U_1,t) - Y, \quad \forall t\in (U_1,U_2]. 
\end{equation}
By the definition of $\eva$ (Eq.~\eqref{eq:E2def}), Eq.~\eqref{eq:QtU1U2}, and the fact that $Q(U_1)\geq 0$, we have that
\begin{equation}
\label{eq:QtU1U2b}
\pb\left(Q(t)\geq [\lambda-(1-p)-\epsilon]t-\zeta-Y \,\big|\, \eva\cap \evb\right)=1, \quad \forall t\in (U_1,U_2]. 
\end{equation}
Let $V$ be the {\it last} time in $(U_1, U_2]$ when the queue length becomes less than $2q_\lambda$, with
\begin{equation}
V = \sup \left\{ t \in [0, B): Q(U_1+t)\leq 2\qlam \right\}, \quad\mbox{if $\inf_{t\in [0, B)}Q(U_1+t)\leq 2\qlam$,}
\label{eq:Vdef}
\end{equation}
and $V=0$, otherwise.  Applying the definition of $V$ in the context of Eq.~\eqref{eq:QtU1U2b} yields that
\begin{equation}
\label{eq:Vlb}
\pb\left(V  \left. \leq \frac{1}{\lambda-(1-p)-\epsilon}\left(2\qlam+Y+ \zeta +1\right) \,\right|\, \eva\cap \evb \right) =1.
\end{equation}
Recall from Proposition \ref{prop:nodel} that, conditional on $\eva\cap\evb$ and assuming $\wlam \asleq c_l\ln \frac{1}{1-\lambda}$, $Y$ must be sub-linear in $B=k\wlam$. In particular, by Eq.~\eqref{eq:Yis0}, we have that, for all $\tau>0$, 
\begin{equation}
\label{eq:Ysmall}
\lim_{\lambda \to 1}\pb\left(Y \leq \tau k \wlam \, \big| \, \eva\cap \evb \right) = 1.
\end{equation}
Combining Eqs.~\eqref{eq:Vlb} and \eqref{eq:Ysmall}, and the fact that $\wlam\to \infty$ as $\lambda \to 1$, we have that, there exists $\upsilon>0$, such that for all $\tau>0$,
\begin{equation}
\label{eq:Vsmall}
\pb\left(V \leq \upsilon\qlam+ \tau k \wlam\, \big| \, \eva\cap \evb \right) = 1-\delta(\lambda), \quad \forall \lambda \in (1-p,1), 
\end{equation}
where $\delta(\cdot)$ is a function with  $\lim_{x \to 1}\delta(x)=0$.  In other words,  conditional on $\eva\cap \evb$, $Q(t)$ will reach the level of $2\qlam$ soon after $U_1$, with high probability. Using the fact that $V\leq U_2$, Eq.~\eqref{eq:Vsmall} further implies that
\begin{equation}
\E\left(V \, \big| \, \eva\cap \evb \right)  \leq  (\upsilon\qlam+ \tau k \wlam)(1-\delta(\lambda))+ U_2\delta(\lambda) \leq \upsilon\qlam+ \tau k\wlam+ U_2\delta(\lambda)
\label{eq:Vsmall2} 
\end{equation}
Translating this into the value of $\E(V)$, we have that
\begin{align}
\limsup_{\lambda \to 1}\frac{\E(V) }{U_2} \leq & \limsup_{\lambda \to 1}\frac{1}{U_2} \big(\E(V\,|\,\eva\cap \evb) \pb(\eva\cap \evb)+{U_2}(1-\pb(\eva\cap \evb))\big)\nnb \\
\stackrel{(a)}{\leq} &\limsup_{\lambda \to 1} \, \left[ 1-\pb(\eva\cap \evb) + \frac{\pb(\eva\cap \evb)}{U_2}\left(\upsilon \qlam+\tau k \wlam + U_2\delta(\lambda)\right) \right] \nnb \\
\stackrel{(b)}{=}& \limsup_{\lambda \to 1} \, \left(  1-\pb(\eva\cap \evb) +\frac{k\wlam}{U_2} \tau \pb(\eva\cap \evb) \right) \nnb\\
\stackrel{(c)}{=}& \limsup_{\lambda \to 1} \, \left(  1-\pb(\eva\cap \evb) +\frac{k}{k+1} \tau \pb(\eva\cap \evb) \right)  \nnb\\
\leq & \tau +  \limsup_{\lambda \to 1} \, \left(  1-\pb(\eva\cap \evb) \right),
\label{eq:EVlim}
\end{align}
where step $(a)$ follows from Eq.~\eqref{eq:Vsmall2},  $(b)$ from the assumptions that $\qlam \ll \wlam$ and  $\lim_{\lambda \to 1}\delta(\lambda)=0$, and $(c)$ from the fact that $U_2=B+\wlam=(k+1)\wlam$. We now connect the behavior of $\E(V)$ to that of $\E(L(0)) = \pb(Q(0)\leq  2\qlam)$, as follows. Fixing any $\lambda \in (1-p, 1)$, we have that
\begin{align}
\E\left(L(0)\right)  \stackrel{(a)}{=} & \E\left(\frac{1}{U_2}\int_{t=0}^{U_2}L(t) dt\right) \nnb \\
\stackrel{(b)}{=}  &\E\left(\frac{1}{U_2}\int_{t=0}^{U1+V}L(t) dt\right) \nnb \\
\stackrel{(c)}{\leq} &   \E\left(\frac{U_1+V}{U_2} \right) \nnb \\
=& \frac{U_1+\E(V)}{U_2}.
 \label{eq:EL0}
\end{align}
where step $(a)$ follows from the stationarity of the process $L(\cdot)$, which in turn follows from the stationarity of $Q(\cdot)$. Step $(b)$ follows from the fact that $L(t)=0$, for all $ t\in [U_1+V,\, U_2]$, which is a consequence of the definition of $V$ in Eq.~\eqref{eq:Vdef}. Step $(c)$ is based on the fact that $L(t)\leq 1$, a.s. By Eq.~\eqref{eq:EL0}, we have that
\begin{align}
\limsup_{\lambda \to 1}\E\left(L(0)\right) \leq & \limsup_{\lambda \to 1} \frac{U_1+\E(V)}{U_2} \nnb\\
\stackrel{(a)}{=} &\frac{\wlam}{(k+1)\wlam} + \limsup_{\lambda \to 1} \frac{\E(V)}{U_2} \nnb \\
\stackrel{(b)}{\leq} &\frac{1}{k}+ \tau +\limsup_{\lambda \to 1} (1-\pb(\eva\cap \evb)),
\end{align}
where steps $(a)$ and $(b)$ follow from the fact that $B=k\wlam$, and Eq.~\eqref{eq:EVlim}, respectively.

By Claim $3$ of Lemma \ref{lem:eveprop}, and Claim $1$ of Lemma \ref{lem:eveIndp}, we have that
\begin{equation}
 \liminf_{\lambda \to 1}\pb\left(\eva\cap \evb\right) = \liminf_{\lambda \to 1}\pb\left(\eva\right)\pb\left(\evb\right)\geq \frac{5}{6}\theta, 
\end{equation}
where $\theta$ is given in Eq.~\eqref{eq:E2isbig}. Set $\tau=k=24$, and let $\zeta$ be sufficiently large so that $\theta\geq 10/9$. We have that
\begin{equation}
	\limsup_{\lambda \to 1}(1-\pb\left(\eva\cap \evb\right))\leq 1-\frac{5}{6}\cdot \frac{9}{10} = 1/4.
\end{equation}
From Eq.~\eqref{eq:EL0}, we have that
\begin{equation}
\limsup_{\lambda \to 1} \E\left(L(0)\right) \leq \frac{1}{k}+ \tau + (1-\pb(\eva\cap \evb)) \leq \frac{1}{24}+\frac{1}{24}+\frac{1}{4}=\frac{1}{3}, 
\end{equation} 
which completes the proof of Proposition \ref{prop:D1lowb}. 
\qed

\subsection{Proof of Theorem \ref{thm:lowerbound}} 
\label{sec:finishingProof}
We now complete the proof of Theorem \ref{thm:lowerbound}. Assuming the validity of Hypothesis \ref{ass:main},  Proposition \ref{prop:D1lowb} asserts that there exists $\cl>0$, so that if $\wlam \asleq \cl\ln\frac{1}{1-\lambda}$ as $\lambda \to 1$, we must have that $\limsup_{\lambda \to 1}\E(L(0))\leq 1/3$ under any sequence of optimal stationary policies. However, this contradicts the requirement that $\E(L(0))\geq {1}/{2}$, given in Eq.~\eqref{eq:visitFreqBound}, which holds independently of the validity of Hypothesis \ref{ass:main}. Therefore, we conclude that Hypothesis \ref{ass:main} must be invalid. 

The invalidity of Hypothesis \ref{ass:main} establishes the lower bound in Eq.~\eqref{eq:qoplb}, as follows.  The negation of the statement of Hypothesis \ref{ass:main} directly implies that there exists $c_l>0$, so that if $\wlam \asleq \cl \ln \frac{1}{1-\lambda}$, as $\lambda \to 1$, then, for any sequence $\{\lambda_n\}$ in $(1-p, 1)$, with $\lim_{n\to \infty}\lambda_n = 1$, we have that
\begin{equation}
\limsup_{n\to \infty} \frac{\mathcal{Q}^*(\lambda_n, W_{\lambda_n})}{\ln \frac{1}{1-\lambda_n}} >0. 
\label{eq:limsup1}
\end{equation}
We can further strengthen Eq.~\eqref{eq:limsup1}, and claim that, for any such sequence, we also have that 
\begin{equation}
\liminf_{n\to \infty} \frac{\mathcal{Q}^*(\lambda_n, W_{\lambda_n})}{\ln \frac{1}{1-\lambda_n}} >0. 
\label{eq:limsup2}
\end{equation}
To show Eq.~\eqref{eq:limsup2}, suppose, for the sake of contradiction, that $\liminf_{n\to \infty} \frac{\mathcal{Q}^*(\lambda_n, W_{\lambda_n})}{\ln \frac{1}{1-\lambda_n}} =0$, for some sequence $\{\lambda_n\}$. This implies that $\{\lambda_n\}$ admits a subsequence, $\{\lambda_{n_k}\}$, such that $\limsup_{k\to \infty} \frac{\mathcal{Q}^*(\lambda_{n_k}, W_{\lambda_n})}{\ln \frac{1}{1-\lambda_{n_k}}} =0$. The existence of the sequence $\{\lambda_{n_k}\}$ contradicts Eq.~\eqref{eq:limsup1}. This proves Eq.~\eqref{eq:limsup2}, which in turn establishes the lower bound in Eq.~\eqref{eq:qoplb}, i.e., that  if $\wlam \asleq \cl \ln \frac{1}{1-\lambda}$, as $\lambda \to 1$, then
\begin{equation}
\qop \ggeq \ln\left(\frac{1}{1-\lambda}\right), \quad \mbox{as $\lambda \to 1$.}
\end{equation}
Finally, we show that the lower bound in Eq.~\eqref{eq:qoplb} is achievable, i.e., that
\begin{equation}
\qop  \lleq  \ln\left(\frac{1}{1-\lambda}\right), \quad \mbox{as $\lambda \to 1$}, 
\label{eq:qupper}
\end{equation}
when $\wlam \asleq \cl \ln\frac{1}{1-\lambda}$. To this end, we invoke Theorem 7 in \cite{SSX12}, which shows that a deterministic queue-length-based diversion policy can achieve the scaling of Eq.~\eqref{eq:qupper}, even when $\wlam=0$.\footnote{As is described in  \cite{SSX12}, the scaling in Eq.~\eqref{eq:qupper} can be achieved by the following simple threshold policy: divert the arrival if and only if the current queue length is equal to a threshold value $x$, where $x$ is set to be the smallest value  such that the resultant rate of diversion is no more than $p$. Since the queue length process under this policy is simply a birth-death process truncated at state $x$, it is easy to verify, via a direct calculation of steady-state probabilities of $Q(t)$, that $\qlam \sim \ln\frac{1}{1-\lambda}$, as $\lambda \to 1$.} This completes the proof of Theorem \ref{thm:lowerbound}.  \qed

\section{Conclusions and Future Work}
\label{sec:conclusion}

In the context of a class of queueing admission control problems, we showed that a non-trivial amount of future information is {necessary} in order to achieve superior heavy-traffic delay performance compared to an online policy. Theorem \ref{thm:lowerbound} also resolves a conjecture posed by \cite{SSX12}. Our proof  exploited certain excursion properties of a transient random walk, which allowed us to connect a policy's diversion decisions to subsequent system idling. 

There are several interesting avenues of future research. First, in light of Theorem \ref{thm:lowerbound} and the results of \cite{SSX12} (Eq.~\eqref{eq:futUpper}), an immediate question is whether the constants $\ch$ and $\cl$ in the scaling of $\wlam$ {coincide}. The granularity of our proof technique does not appear to be sufficient to answer this question, which likely demands a finer analysis. 

Because our proof relies mostly on the {macroscopic} properties of the input sample paths, the techniques and resultant insights in this paper seem to be fairly robust and can potentially be generalized to derive lower bounds on the necessary amount of future information for other resource allocation problems. For example, one generalization could be for a setting where the arrival and service token processes are non-Poisson (e.g., renewal or phase-type processes). In this case, we expect similar arguments to work when the process, $S(0, \cdot)$, admits similar excursion properties as in the case of Poisson processes, and does not exhibit substantial long-range correlations (for otherwise, one could potentially obtain more future information by looking into the history of past inputs). Another possibility would be to consider systems with multiple queues, in which case the relevant excursion properties of the input processes would likely be connected to those of random walks in higher dimensions. Yet another variation would be to relax the hard diversion rate constraint, and consider instead the scenario where the system manager is interested in minimizing some combined cost as a function of the delay and diversion rate. However, depending on the cost function, one may need to adjust the performance metric or regime of interest, since the system may not ever have to become critically loaded, simply because the cost structure would encourage a higher rate of diversion as the system load increases. 

Finally, at a higher level, while our result focuses on the \emph{quantity} of future information, measured by the length of a lookahead window, there is another important dimension of \emph{quality}. For instance, the observed future input may differ from the actual realizations due to prediction noise, or alternatively, only distributional information of future input is available. Neither our results, nor those of \cite{SSX12},  deal with the impact of prediction noise, and \cite{XuChan14} considers only a specific noise model induced by random no-shows. A rigorous understanding of the impact of prediction accuracy in the context of dynamic resource allocation problems could be a promising direction for future research.

\section{Acknowledgment}

The research reported in this paper was supported in part by NSF grant CMMI-1234062 and a Claude E.~Shannon Research Assistantship. The author would like to thank the anonymous reviewers for their detailed and insightful comments. 

\bibliographystyle{ormsv080}
\bibliography{LB_Bib} 

\begin{thebibliography}{15}
\expandafter\ifx\csname natexlab\endcsname\relax\def\natexlab#1{#1}\fi
\expandafter\ifx\csname url\endcsname\relax
  \def\url#1{{\tt #1}}\fi
\expandafter\ifx\csname urlprefix\endcsname\relax\def\urlprefix{URL }\fi
\expandafter\ifx\csname urlstyle\endcsname\relax
  \expandafter\ifx\csname doi\endcsname\relax
  \def\doi#1{doi:\discretionary{}{}{}#1}\fi \else
  \expandafter\ifx\csname doi\endcsname\relax
  \def\doi{doi:\discretionary{}{}{}\begingroup \urlstyle{rm}\Url}\fi \fi

\bibitem[{Brown et~al.(2010)Brown, Smith, and Sun}]{Brown10}
Brown, D.~B., J.~E. Smith, P.~Sun. 2010.
\newblock Information relaxations and duality in stochastic dynamic programs.
\newblock {\it Operations Research\/} {\bf 58}(4) 785--801.

\bibitem[{Desai et~al.(2012)Desai, Farias, and Moallemi}]{Farias12}
Desai, V.~V., V.~F. Farias, C.~C. Moallemi. 2012.
\newblock Pathwise optimization for optimal stopping problems.
\newblock {\it Management Science\/} {\bf 58}(12) 2292--2308.

\bibitem[{Fisher and Raman(1996)}]{FR96}
Fisher, M., A.~Raman. 1996.
\newblock Reducing the cost of demand uncertainty through accurate response to
  early sales.
\newblock {\it Operations Research\/} {\bf 44}(1) 87--99.

\bibitem[{Gonzlez-Hern\'{a}ndez and Villarreal(2011)}]{GV11}
Gonzlez-Hern\'{a}ndez, J., C.~E. Villarreal. 2011.
\newblock Optimal policies for constrained average-cost {M}arkov decision
  processes.
\newblock {\it TOP\/} {\bf 19}(1) 107--120.

\bibitem[{Hern\'{a}ndez-Lerma et~al.(2003)Hern\'{a}ndez-Lerma,
  Gonzalez-Hern\'{a}ndez, and L\'{o}pez-Martinez}]{HGL03}
Hern\'{a}ndez-Lerma, O., J.~Gonzalez-Hern\'{a}ndez, R.~R. L\'{o}pez-Martinez.
  2003.
\newblock Constrained average cost {M}arkov control processes in {B}
  spaces.
\newblock {\it SIAM Journal on Control and Optimization\/} {\bf 42}(2)
  442--468.

\bibitem[{Kim and Horowitz(2002)}]{KH02}
Kim, S.~C., I.~Horowitz. 2002.
\newblock Scheduling hospital services: The efficacy of elective surgery
  quotas.
\newblock {\it Omega\/} {\bf 30} 335--346.

\bibitem[{Rogers(2007)}]{Rogers07}
Rogers, L. C.~G. 2007.
\newblock Pathwise stochastic optimal control.
\newblock {\it SIAM J. Control Optim.\/} {\bf 46}(3) 1116--1132.

\bibitem[{Spencer et~al.(2014)Spencer, Sudan, and Xu}]{SSX12}
Spencer, J., M.~Sudan, K.~Xu. 2014.
\newblock Queuing with future information.
\newblock {\it Annals of Applied Probability\/} {\bf 24}(5) 2091--2142.

\bibitem[{Stidham(1985)}]{Sti85}
Stidham, S.Jr. 1985.
\newblock Optimal control of admission to a queueing system.
\newblock {\it IEEE Trans.\ Automatic Control\/} {\bf 30}(8) 705--713.

\bibitem[{Sun et~al.(2009)Sun, Heng, Seow, and Seow}]{sun_09}
Sun, Yan, Bee~H Heng, Yian~T Seow, Eillyne Seow. 2009.
\newblock Forecasting daily attendances at an emergency department to aid
  resource planning.
\newblock {\it BMC emergency medicine\/} {\bf 9}(1) 1.

\bibitem[{Tsitsiklis and Xu(2012)}]{TX12}
Tsitsiklis, J.~N., K.~Xu. 2012.
\newblock On the power of (even a little) resource pooling.
\newblock {\it Stochastic Systems\/} {\bf 2}(1) 1--66.

\bibitem[{Wargon et~al.(2009)Wargon, Guidet, Hoang, and Hejblum}]{wardon_emj09}
Wargon, M., B.~Guidet, T.D. Hoang, G.~Hejblum. 2009.
\newblock A systematic review of models for forecasting the number of emergency
  department visits.
\newblock {\it Emergency Medicine Journal\/} {\bf 26}(6) 395--399.

\bibitem[{Xu(2014)}]{xu2014power}
Xu, K. 2014.
\newblock On the power of (even a little) flexibility in dynamic resource
  allocation.
\newblock Ph.D. thesis, Massachusetts Institute of Technology.

\bibitem[{Xu and Chan(2014)}]{XuChan14}
Xu, K., C.~W. Chan. 2014.
\newblock Using future information to reduce waiting times in the emergency
  department via diversion.
\newblock {\it Manuscript\/} .

\bibitem[{Yechiali(1971)}]{Yech71}
Yechiali, U. 1971.
\newblock On optimal balking rules and toll charges in the {$GI/M/1$} queuing
  process.
\newblock {\it Operations Research\/} {\bf 19}(2) 349--370.

\end{thebibliography}

\begin{APPENDIX}{}

\section{Additional Proofs}

\subsection{Proof of Lemma \ref{lem:stationaryStabel}}
\label{app:lem:stationaryStable}

\def\calX{\mathcal{X}}
\def\calY{\mathcal{Y}}

\bproof  We will formulate our admission control problem as a discrete-time Markov decision process (MDP), and  invoke existing results to verify the existence of a stable stationary optimal policy. Recall that state of the system at the $n$th step is $X_n=(Q(Z_n-), F_n)$, where $F_n$ was defined in Eq.~\eqref{eq:defFn}. Define $\mathcal{X}$ as the set
\begin{equation}
\mathcal{X} = \zp\times [-1, w]^{\N}. 
\label{eq:XspaceDef}
\end{equation}
Note that $X_n$ can be represented as an element in $\mathcal{X}$ for all $n$, because $Q(Z_n-)$ is the queue length just before the $n$th event and hence belongs to $\zp$, and each coordinate of $F_n$, which either represents the type of an event or an inter-arrival time upper-bounded by $\wlam$, lies in the interval $[-1, \wlam]$. The following topological properties of $\calX$ are useful, whose proof is given in Appendix \ref{app:lem:Xtopo}. 

\begin{lemma}
\label{lem:Xtopo} The following holds. 
\begin{enumerate}
\item $\calX$ is Polish, i.e., it is complete and separable. 
\item Under an appropriate metric, the set $\{x\in \calX : x_1 \leq a\}$ is compact for all $a \in \rp$. 
\end{enumerate}
\end{lemma}
The MDP associated with our admission control problem is defined as follows: 
\begin{enumerate}
\item The state space is $\mathcal{X}$, defined in Eq.~\eqref{eq:XspaceDef}. 
\item The action space, $\mathcal{L}$, is the closed interval $[0,1]$, and the action at step $n$, $l_n \in \mathcal{L}$ ,  specifies the probability of diversion, i.e., $l_n = \pb(H(n)=1)$. Denote by $\mathcal{L}(X)$ the set of allowable actions when the system is in state $X$. Then $\mathcal{L}(X_n)=[0,1]$ if $A_n(0)=1$, which corresponds to the $n$th event being an arrival, and $\mathcal{L}(X_n)=\{0\}$ if $S_n(0)=1$, which corresponds to the $n$th event being the generation of a service token.  
\item The stochastic kernel is the one associated with the Poisson arrival and service token processes, as well as the queueing and diversion dynamics. 
\item The $n$th step is associated with a \emph{penalty}, $f(X_n, l_n)$, which is equal to the queue length, $Q(Z_n-)$. It also incurs a \emph{cost}, $c(X_n, l_n)$, which is equal to the probability of diversion, $l_n$. 
\item The objective is to minimize the {time-average penalty}, defined in Eq.~\eqref{eq:avgQueue}, subject to a constraint on the {time-average cost}, defined in Eq.~\eqref{eq:avgCost}. 
\end{enumerate}
Theorem 3.2 and Lemma 3.5 of \cite{HGL03} show that  an MDP of this kind admits a stable stationary optimal policy, provided that a set of conditions are satisfied, which are given in Section 2 and Assumption 3.1 of \cite{HGL03}. These conditions are met by our MDP, and we highlight a few among them: $(1)$ the state space is Polish (by the first claim of Lemma \ref{lem:Xtopo}), $(2)$ the set $\{(X, l) \in (\mathcal{X}, \mathcal{L}): f(X_n, l_n) \leq a\}$ is compact for all $a \in \rp$ (by the second claim of Lemma \ref{lem:Xtopo}), $(3)$ $c(X_n, l_n)$, which in our case is simply equal to $l_n$, is non-negative and lower semi-continuous in $l_n$ for every state $X_n \in \mathcal{X}$, and $(4)$ the stochastic kernel satisfies a certain weak continuity condition, which essentially requires the distribution of $X_n$ not vary abruptly as a function of the state-action pair $(X_n, l_n)$, and this continuity condition can be verified by using the definitions of  Poisson processes and the associated queueing dynamics. This completes the proof of Lemma \ref{lem:stationaryStabel}. 
\qed

\subsection{Proof of Lemma \ref{lem:eveprop}}
\label{app:lem:eveprop}

\normalsize

\bproof
Recall from Eq.~\eqref{eq:Sdef} that $S(s,t)$ is defined as the difference between the numbers of arrivals and service tokens in $(s,t]$. Since the arrival and service tokens processes are independent Poisson processes with rate $\lambda$ and $1-p$, respectively, it is not difficult to verify that 
\begin{equation}
S(s,t) \stackrel{d}{=} \sum_{n = 1}^{N_{s,t}} X_n,
\label{eq:Sdef2}
\end{equation}
where $N_{s,t}$ is a Poisson random variable with mean $(\lambda+1-p)(t-s)$, which corresponds to the total number of events in $(t,s]$, and the $X_n$s are i.i.d., with 
\begin{equation}
X_1 = \left\{
\begin{array}{ll}
1,& \quad \mbox{w.p. } \frac{\lambda}{\lambda+1-p}, \\
-1,& \quad  \mbox{otherwise},
\end{array}
\right.
\end{equation}
By Eq.~\eqref{eq:Sdef2}, and the fact that $\lim_{B\to \infty} \frac{N_{s,s+B} }{B}= \lambda+1-p$ almost surely, Claim $1$ follows from a variation of the standard Functional Law of Large Numbers (FLLN) for the sum of bounded i.i.d.~random variables. Claim $3$ follows from the Weak Law of Large Numbers applied to the sum of i.i.d.\ Poisson random variables, and our assumption that $W_\lambda \to \infty$ as $\lambda \to 1$ (Eq.~\eqref{eq:Wlamassump}).  Finally, Claim $2$ follows from the Markov's inequality, in the same way as in Eq.~\eqref{eq:visitFreqBound}, by noting that $\E(Q(0))=\qlam$ under an optimal stationary policy.  \qed

\subsection{Proof of Lemma \ref{lem:Eve}}
\label{app:lem:PoissonLDP}

\bproof
Based on the stationarity of $\arr$ and $\ser$, and the assumption that $B=k\wlam$ and $\qlam\ll \wlam$, it suffices for us to show,  that for any $a,b>0$, there exists $\gamma>0$, such that	
\begin{equation}
\pb\left(S(0,a \wlam)\leq - b \wlam \right) \ggeq \exp(-\gamma \wlam), \quad \mbox{as $\lambda \to 1$}.
\label{eq:eveEqv}
\end{equation}

By definition, the distribution of $S(0,t)$ can be written as 
\begin{equation}
S(0,t) \stackrel{d}{=}  A_{\lambda t}- D_{(1-p)t},
\end{equation}
where $A_{\lambda t}$ and $ D_{(1-p)t}$ are independent Poisson random variables with mean $\lambda t$ and $(1-p)t$, respectively. The following lemma follows from the standard large-deviation principles of Poisson random variables, and its proof is omitted. 
\begin{lemma}
\label{lem:PoissonLDP}
Let $D_x$ be a Poisson random variable with mean $x$. Then, for all $c_1>0$, there exists $c_2>0$, such that 
\begin{equation}
\pb\left(D_x \geq c_1 x\right) \ggeq  \exp(-c_2 x), \quad \mbox{as $x \to \infty$}.
\end{equation}
\end{lemma}

Combining Lemma \ref{lem:PoissonLDP} and the fact that $\wlam\to \infty$ as $\lambda \to 1$, we have that there exists $\gamma>0$, such that
\begin{equation}
\pb\left(D_{(1-p)a \wlam} \geq (b+2a)\wlam\right)\ggeq \exp(-\gamma \wlam)
\label{eq:DisBig}
\end{equation}
as $\lambda \to 1$. We have that
\begin{align}
& \pb\left(S(0,a \wlam)\leq - b \wlam \right) \nnb\\
\geq & \pb\left(\left\{A_{\lambda a\wlam} < 2a\wlam\right\}\cap \left\{D_{(1-p)a\wlam}\geq (b+2a)\wlam\right\}\right) \nnb\\
\stackrel{(a)}{=}  &\pb\left({A_{\lambda a\wlam} < 2a\wlam}\right)\pb\left( {D_{(1-p)a\wlam}\geq (b+2a)\wlam}\right)\nnb\\
\stackrel{(b)}{\geq}  &\pb\left({A_{\lambda a\wlam} < 2 \lambda a\wlam}\right)\pb\left( {D_{(1-p)a\wlam}\geq (b+2a)\wlam}\right)\nnb\\
\stackrel{(c)}{\geq}  &\frac{1}{2}\pb\left( {D_{(1-p)a\wlam}\geq (b+2a)\wlam}\right)\nnb\\
\stackrel{(d)}{\ggeq}  &\exp(-\gamma\wlam),
\end{align}
as $\lambda \to 1$, where step $(a)$ follows from the independence between $A_{\lambda a\wlam}$ and $D_{(1-p)a\wlam}$, $(b)$ from the fact that $\lambda<1$, $(c)$ from the Markov's inequality, and $(d)$ from Eq.~\eqref{eq:DisBig}. This proves Eq.~\eqref{lem:Eve}, and hence Lemma \ref{lem:Eve}.  \qed

\subsection{Proof of Lemma \ref{lem:eveIndp}}
\label{app:lem:eveIndp}

\bproof For Claim 1, observe that each of the event concerns only the behavior of the arrival and service token processes over an interval, and that these intervals are disjoint from each other. Claim $1$ follows by noting that both $\arr$ and $\ser$ are Poisson processes and hence memoryless. For Claim $2$, because the policy has access to a lookahead window of length $\wlam$, the queue length at time $t$ is hence $\calf_{t+\wlam}$ measurable, where $\calf$ is the natural filtration induced by the input processes. The claim follows again from the memoryless property of Poisson processes. Claim $3$ follows from the same arguments as for Claim $2$. 
\qed

\subsection{Proof of Lemma \ref{lem:EJLB}}
\label{app:lem:EJLB}

\bproof Consider the sequence of optimal stationary policies, $\{\pi_\lambda\}$. Let $\phi$ be defined as in Eq.~\eqref{eq:E6def}. Fix $\phi >0$, and let 
\begin{equation}
	K = U_3 +\phi \wlam \stackrel{(a)}{=} (k+\phi+2)\wlam, 
\end{equation}
where step $(a)$ follows from the fact that $U_3 = B+2\wlam$ and  $B=k\wlam$. 
The main idea for the proof is based on the following observation: conditional on $\cap_{i=1}^5 \mathcal{E}_i$, the queue length process, $Q(t)$,  would have reached zero before time $K$, even if \emph{no} diversion had been made in $(0, K]$ (illustrated in Figure \ref{fig:samplePath}). Therefore, each diversion made in $(U_1, U_2]$ will necessarily lead to a \emph{waste service token} in $(0,K]$, and hence
\begin{align}
\label{eq:JYcouple}
	& \pb\left( \jc(K) \geq \tau B  \bsep \cap_{i=1}^5 \mathcal{E}_i \right) \geq \pb\left(  Y\geq \tau B \,\big|\, \cap_{i=1}^5 \mathcal{E}_i \right). 
	\end{align}
We next give a lower bound on the above probability, as follows: 
\begin{align}
   &\pb\left( \jc(K) \geq \tau B\bsep \eva \cap \evb  \right)  \nnb \\
   \geq& \pb\left(\jc(K)\geq \tau B , \cap_{i=3}^5\mathcal{E}_i \bsep \eva \cap \evb\right) \nnb\\
   =&\pb\left(\jc(K)\geq \tau B \bsep \cap_{i=1}^5\mathcal{E}_i \right) \pb\left(\cap_{i=3}^5 \mathcal{E}_i \bsep \eva \cap \evb \right) \nnb\\
   \sk{a}{\geq}&\pb\left(Y \geq \tau B \bsep \cap_{i=1}^5\mathcal{E}_i \right) \pb\left(\cap_{i=3}^5 \mathcal{E}_i \bsep \eva \cap \evb \right) \nnb\\
   = &\pb\left(Y \geq \tau B , \cap_{i=3}^5\mathcal{E}_i \bsep \eva \cap \evb\right) \nnb\\
   \sk{b}{=}& \pb\left(\eve\right)\pb\left(Y\geq \tau B, \evc \cap \evd \bsep \eva \cap \evb \right)\nnb\\
   \geq & \pb\left(\eve\right)\left(\pb\left(Y\geq \tau B\bsep \eva \cap \evb \right) + \pb\left(\evc \bsep \eva\cap \evb\right)+\pb\left(\evd \bsep \eva \cap \evb \right)-2\right)\nnb\\
   \geq & \pb\left(\eve\right)\left(\pb\left(Y\geq \tau B\bsep \eva \cap \evb \right) + \pb\left(\evc \bsep \eva\cap \evb\right)+\frac{\pb\left(\evd\right)+\pb\left(\eva\cap \evb\right)-1}{\pb\left(\eva\cap \evb\right)}-2\right)\nnb\\
   \sk{c}{=} & \pb\left(\eve\right)\left(\pb\left(Y\geq \tau B\bsep \eva \cap \evb \right) + \pb\left(\evc\right)+\frac{\pb\left(\evd\right)-1}{\pb\left(\eva \cap \evb \right)}-1\right)
   \label{eq:Jkpb1}
\end{align}
where step $(a)$ follows from Eq.~\eqref{eq:JYcouple}, and $(b)$ and $(c)$ from the independence between $\eve$ and $\eva\cap \evb$, and between $\evc$ and $\eva\cap \evb$, respectively (Lemma \ref{lem:eveIndp}). We have also used the inequality that $\pb\left(A\cap B\right)\geq \pb\left(A\right)+\pb(B)-1$, for any events $A$ and $B$. 

By Claim $3$ of Lemma \ref{lem:eveprop},  we have that
\begin{equation}
\label{eq:ev34conv}
\lim_{\lambda \to 1}\pb\left(\evc\right)=\lim_{\lambda \to 1}\pb\left(\evd\right)=1. 
\end{equation}
Combining the assumption (Eq.~\eqref{eq:Ycontra1})
\begin{equation}
\liminf_{\lambda \to 1}\pb\left(Y\geq \tau B \,\big|\, \eva \cap \evb  \right)=q>0 
\end{equation}
with Eqs.~\eqref{eq:Jkpb1} and \eqref{eq:ev34conv}, we have that there exists $\tilde{\lambda}\in(0,1)$, such that
\begin{align}
\pb\left( \jc(K) \geq \tau B\bsep \eva \cap \evb  \right) \geq& \pb\left(\eve \right)\pb\left( Y \geq \tau B \bsep \eva \cap \evb \right) \nnb\\
\geq & \pb\left(\eve\right)q/2,
\label{eq:JkLB2}
\end{align}
for all $\lambda \in \left(\tilde{\lambda}, 1\right)$.  
We have that 
\begin{align}
	\E\left( \jc(K)\right) 
	\geq&  \tau B \cdot \pb\left( \jc(K) \geq \tau B \right) \nnb \\
	\geq & \tau B \cdot \pb\left( \jc(K)\geq \tau B, \eva\cap \evb\right) \nnb \\ 
	=&\tau B \cdot \pb\left( \jc(K)\geq \tau B \,\big|\, \eva\cap \evb\right) \cdot \pb\left( \eva\cap \evb\right)\nnb \\
	\stackrel{(a)}{\ggeq}& B  \pb\left(\eve\right) \pb\left( \eva\cap \evb\right) \nnb \\
	\stackrel{(b)}{\ggeq}&  B  \pb\left(\eve\right) \nnb\\
	\sk{c}{\ggeq}& B \exp\left(-\gamma\wlam\right) ,
	\label{eq:tailLower}
\end{align}
for some $\gamma >0$, as $\lambda \to 1$, where step $(a)$ follows from Eq.~\eqref{eq:JkLB2}, $(b)$ from Claims $1$ and $2$ of  Lemma \ref{lem:eveprop} and the independence of the events $\eva$ and $\evb$ (Claim $1$ of Lemma \ref{lem:eveIndp}), i.e., that
\begin{equation}
\label{eq:Pevab}
	\pb\left(\eva \cap \evb \right) = \pb\left(\eva\right)\pb\left(\evb\right) \geq \frac{5}{6}\theta, 
\end{equation}
and $(c)$ from Lemma \ref{lem:Eve}. This proves Lemma \ref{lem:EJLB}, by setting $a = k+\phi+2$.  $\qed$

\subsection{Proof of Lemma \ref{lem:Xtopo}}
\label{app:lem:Xtopo}

\bproof  Let $\calX_0 = [-1, \wlam]^\N$. We will show that $\calX_0$ is compact under the metric $\|x-y\|_g = \sum_{i=1}^\infty 2^{-i}|x_i-y_i|$. If this is true, it is not difficult to show that, for any $a \in \rp$, the set $\{x\in \calX : x_1 \leq a\} = \{ 0,\ldots, \lfloor a \rfloor \}\times \calX_0$ is also compact under $\|\cdot\|_g$, and our second claim follows. Note that a compact metrizable space is Polish, and it is easy to show that $\zp$ is Polish under the $l_1$ norm. Our first claim thus also follows from the compactness of $\calX_0$, by observing that the product of two Polish spaces remains Polish. 

We now show the compactness of $\calX_0$. It suffices to show that any sequence in $\calX_0$,  $\{x^{i}\}_{i\in \N}$, admits a sub-sequence that converges to a point in $\calX_0$. We will construct such a limiting point coordinate-by-coordinate, as follows. Because $x^{i}_1$ is an element of the compact interval $[-1,\wlam]$ for all $i\in \N$, there exists $y_1\in [-1,\wlam]$ and an increasing sequence, $\{i^{1,j}\}_{j\in \N} \subset \N$, such that $\lim_{j\to \infty} x^{i^{1,j}}_1 = y_1$. We now apply the same reasoning for progressively larger values of $k$: there exist  $y_k \in [-1,\wlam]$ and  $\{i^{k,j}\}_{j\in \N}$ for $k=2,3,\ldots$, such that, for every $k\geq 2$, $\{i^{k,j}\}_{j\in \N}$ is a sub-sequence of $\{i^{k-1,j}\}_{j\in \N}$, and 
\begin{equation}
 \lim_{j\to \infty}x^{i^{k,j}}_k= y_k.
 \label{eq:ylim1}
\end{equation} 
Fix $k\geq 2$. Because $\{i^{k,j}\}_{j\in \N}$ is a sub-sequence of $\{i^{m,j}\}_{j\in \N}$ for all $m\leq k-1$, Eq.~\eqref{eq:ylim1} further implies that
\begin{equation}
 \lim_{j\to \infty}x^{i^{k,j}}_m= y_m, \quad \forall m\in \{1,\ldots, k\}, 
\end{equation} 
or, equivalently, that
\begin{equation}
\lim_{j\to \infty} \sum_{i=1}^k 2^{-m} \left|x_m^{i^{k,j}} - y_m \right | = 0, \quad \forall k\in \N. 
 \label{eq:yklimit}
\end{equation}

Let $y$ be the element of $\calX_0$ whose coordinates are defined according to the above procedure. We argue that $y$ is the limiting point for some sub-sequence of $\{x^{i}\}_{i=1}^\N$. For every $k \in \N$, there exists $j(k)\in \N$, such that for all $j\geq j(k)$, 
\begin{align}
\left\| y-x^{i^{k, j}} \right \|_g =&  \sum_{m=1}^\infty2^{-m} \left| y_m-x^{i^{k, j}}_m  \right|  \nnb\\
\stackrel{(a)}{\leq} & \left( \sum_{m=1}^k 2^{-m} \left| y_m-x^{i^{k, j}}_m  \right| \right) +(1+\wlam)2^{-(k-2)}  \nnb\\
\stackrel{(b)}{\leq} & \frac{1}{k} +(1+\wlam)2^{-(k-2)}, 
\label{eq:yminusxShrink}
\end{align}
where step $(a)$ follows from the fact that $\left| y_m-x^{i^{k, j}}_m  \right|\leq 2(\wlam+1)$ for all $m \in \N$, and step $(b)$ from Eq.~\eqref{eq:yklimit}. Define
\begin{equation}
n^k = \max \{i^{m,j(m)}: 1\leq m\leq k\}, \quad \forall k\in \N. 
\end{equation}
By Eq.~\eqref{eq:yminusxShrink}, we have that
\begin{equation}
\left\| y-x^{n^k} \right \|_g  \leq \frac{1}{k} +(1+\wlam)2^{-(k-2)}, \quad \forall k\in \N, 
\end{equation}
Therefore, $\{x^{n^k}\}_{k\in \N}$ is a sub-sequence of $\{x^{i}\}_{i \in \N}$, and it converges to $y$ as $k\to \infty$ under the metric $\|\cdot\|_g$. This proves that $\calX_0$ is compact. 
\qed

\end{APPENDIX}

\end{document}